\documentclass[12pt,leqno]{amsart}
\usepackage{amssymb}
\usepackage{longtable}
\usepackage[all]{xy}
\usepackage{amsxtra}

\newtheorem{theorem}{Theorem}
\newtheorem{corollary}[theorem]{Corollary}

\newtheorem{proposition}[theorem]{Proposition}

\theoremstyle{definition}
\newtheorem{defi}[theorem]{Definition}
\newtheorem{example}[theorem]{Example}
\newtheorem{remark}[theorem]{Remark}
\newtheorem*{ack}{Acknowledgement}

\numberwithin{theorem}{section}

\numberwithin{equation}{section}

 \DeclareMathOperator{\aut}{Aut}
 
 \DeclareMathOperator{\spp}{Sp}
 \DeclareMathOperator{\psl}{PSL}

 \DeclareMathOperator{\diag}{diag}

\title{Weyl Groups and Abelian Varieties}
\author[A. Carocca]{Angel Carocca}
\address{Facultad de Matem\'aticas, Pontificia Universidad Cat\'olica
de Chile, Casilla 306, Correo 22, Santiago, Chile}
\email{acarocca@mat.puc.cl}
\author[V. Gonz\'alez-Aguilera]{V\'\i ctor Gonz\'alez-Aguilera}
\address{Departamento de Matem\'aticas, Universidad T\'ecnica
Federico Santa Mar\'\i a, Casilla 110-V, Valpara\'\i so, Chile}
\email{vgonzale@mat.utfsm.cl}
\author[R. E. Rodr\'\i guez]{\newline Rub\'\i \ E. Rodr\'\i guez}
\address{Facultad de Matem\'aticas, Pontificia Universidad Cat\'olica
de Chile, Casilla 306, Correo 22, Santiago, Chile}
\email{rubi@mat.puc.cl}
\thanks{This research was partially supported by the Presidential
Science Chair on Geometry and by FONDECYT Grants \# 1030595,
1040319 and 1011039.}

\begin{document}

\begin{abstract}
Let $G$ be a finite group. For each integral representation $\rho$
of $G$ we consider $\rho-$de\-com\-po\-sable principally polarized
abelian varieties; that is, principally polarized abelian
varieties $(X,H)$ with $\rho(G)-$action, of dimension equal to the
degree of $\rho$, which admit a decomposition of the lattice for
$X$ into two $G-$invariant sublattices isotropic with respect to
$\Im H$, with one of the sublattices $\mathbb{Z}G-$isomorphic to
$\rho$.

We give a construction for $\rho-$de\-com\-po\-sable principally
polarized abelian varieties, and show that each of them is
isomorphic to a product of elliptic curves.

Conversely, if $\rho$ is absolutely irreducible, we show that each
$\rho-$de\-com\-po\-sable p.p.a.v. is (isomorphic to) one of those
constructed above, thereby characterizing them.

In the case of irreducible, reduced root systems, we consider the
natural representation of its associated Weyl group, apply the
preceding general construction, and characterize completely the
associated families of principally polarized abelian varieties,
which correspond to modular curves.
\end{abstract}



\maketitle

\section{Introduction}
\label{sec:intro}
It is known that every finite group acts on some
curve, hence on some Jacobian, and therefore on some principally
polarized abelian variety (\textit{p.p.a.v.}). However, there is
no explicit way of relating the order of the group and the
dimension of the p.p.a.v. on which it acts.

For any integral representation of a finite group $G$ there exists
a family of principally polarized abelian varieties of dimension
\textit{equal to} the degree of the representation, each invariant
under the given action of $G$.

In particular, every root system of rank $n$ has a natural lattice in
${\mathbb R}^n $ associated to it, and hence a natural integral
representation of its associated Weyl group of degree $n$.
Therefore, we may apply the preceding construction to build a
family of p.p.a.v.'s of dimension $n$ admitting the given Weyl
group (more generally, the full automorphism group of the root
system) as subgroup of its group of automorphisms.

In the case of irreducible, reduced root systems, we show that each
p.p.a.v. in the corresponding family is isomorphic (as a complex
torus) to a product of elliptic curves. These results are related
to an invariant theory developed in the work of E.
Looijenga~\cite{loo} where he also deals with abelian varieties
associated to each irreducible root system, all of which admit the
action of the corresponding Weyl group and are isomorphic (as
complex tori) to a product of elliptic curves.

In general, the reducibility or irreducibility of a given p.p.a.v.
(as p.p.a.v.) is not an easy problem. Here we show that, in most
cases, our varieties are irreducible as p.p.a.v.'s.

In addition, in the case of irreducible, reduced root systems, we
characterize the corresponding families in the moduli space of
p.p.a.v.'s $\mathcal{A}_n$, and show that they are uniformized by
explicit modular curves.

In Section \ref{sec:intro} we recall some basic results on
principally polarized abelian varieties and roots systems in order
to fix the notation.

In Section \ref{sec:genconst} we consider any integral
representation $\rho$ of a finite group $G$, and construct
p.p.a.v.'s $(X,H)$, of dimension equal to the degree of $\rho$,
which admit a decomposition of the lattice for $X$ into two
$G-$invariant sublattices isotropic with respect to $\Im H$, with
one of the sublattices $\mathbb{Z}G-$isomorphic to $\rho$. We call
these p.p.a.v.'s $\rho-$de\-com\-po\-sable. We also show that each
$\rho-$decomposable p.p.a.v. is isomorphic to a product of
elliptic curves.

Conversely, if $\rho$ is absolutely irreducible, we show that each
$\rho-$de\-com\-po\-sable p.p.a.v. is (isomorphic to) one of those
constructed above, thereby characterizing them.

In Section \ref{sec:constweyl}, we apply our characterization to
each irreducible, reduced root system and obtain the corresponding
Riemann matrices. We also discuss the relationship between our
varieties and the previous work of Looijenga, and the notion of
dual polarization developed by Birkenhake and Lange \cite{bir}.

Section \ref{sec:irred} describes explicit isomorphisms of our
principally polarized abelian varieties with products of elliptic
curves and discusses their irreducibility as principally polarized
abelian varieties.

Section \ref{sec:modular} shows that our families for irreducible,
reduced root systems are uniformized by modular curves, which are
described explicitly, and complete their connection with the
previous work of Looijenga.

Now we fix the notation as follows. A \textit{complex torus} of
dimension $n$ will be by definition a quotient $X\simeq {\mathbb
C}^n/L$ of the vector space ${\mathbb C}^n$ modulo a lattice $L$
of maximal rank. $X$ is called an \textit{abelian variety} if it
admits a polarization; that is, a positive definite hermitian form
$H$ on ${\mathbb C}^n$ whose imaginary part is integral valued on
$L$. In this case, the lattice $L$ always admits a basis such that
the matrix of $\Im H$ in this basis
is of the form $\left(%
\begin{array}{rr}
  0_n & D \\
  -D & 0_n \\
\end{array}%
\right) $ with a diagonal matrix $D=\diag (d_1,d_2,...,d_n)$ and
positive integers $d_i$ with $d_i/d_{i+1}$. The $n$-tuple $(d_1,
d_2, \ldots , d_n)$ is called the \textit{type} of the
polarization, the integer $d=d_1 d_2 \ldots d_n$ is called the
\textit{degree} of the polarization and a polarization of type
$(1,1,...,1)$ is called a \textit{principal} polarization.

The automorphisms of the complex torus $X$ that respect the
polarization $H$ form a finite group denoted by $\aut_H(X)$. If
$G$ is a subgroup of $\aut_H(X)$, we will say that $(X,H)$ has a
$G-$action. The representation of $\aut_H(X)$ on the lattice of
$X$ will be called the \textit{rational} representation of
$\aut_H(X)$.

A \textit{principally polarized abelian variety} of dimension $n$
(p.p.a.v.) will be a pair $(X,H)$ where $X$ is an $n$-dimensional
complex torus and $H$ is a principal polarization on $X$.
Alternatively we can consider a principally polarized abelian
variety $(X,H)$ as a pair $(X,\mathcal L)$ where $\mathcal{L}$ is
a line bundle on $X$ such that its first Chern class
$c_1(\mathcal{L}) \simeq H$, but in order to give explicit
descriptions we will follow the classical terminology.

Each Riemann matrix $Z$ in the Siegel space  $\mathbb{H}_n$ of
degree $n$ gives a p.p.a.v. of dimension $n$. If we denote by
$\spp(2n, \mathbb{Z})$ the symplectic group, then the complex
analytic space $\mathcal{A}_n = \mathbb{H}_n/ \spp(2n, {\mathbb
Z})$ parametrizes the (isomorphism classes of) p.p.a.v.'s of
complex dimension $n$.

Since there are different kinds of equivalences between Riemann
matrices or abelian varieties, and so as to avoid ambiguities, we
recall the definitions involved in our results.

Two p.p.a.v.'s $(X_1,H_1)$ and $(X_2,H_2)$ are called
\textit{isogenous} if there exist  a surjective homomorphism $h:
X_1 \rightarrow X_2$ (as complex tori) with $0$-dimensional
kernel; they are called \textit{isomorphic} if the kernel of $h$
is trivial. They will be called  \textit{isomorphic as principally
polarized abelian varieties} if there exists an isomorphism $h$ (
as above) which further preserves the given polarizations $H_1$
and $H_2$.

A p.p.a.v. is called \textit{irreducible} if it is not isomorphic
as a \textit{principally polarized abelian variety} to the product
of p.p.a.v.'s of smaller dimensions.

For root systems $R$ we will follow the definitions and notation
in Bourbaki \cite{b}. In particular, $R^{\vee}$ will denote the
inverse root system and  the Weyl group of $R$ will be denoted by
$\mathcal{W}(R)$. We will essentially deal with irreducible and
reduced root systems in the sense of \cite{b}; that is, with the
root systems denoted by $\mathbf{A_n}$, $\mathbf{B_n}$,
$\mathbf{C_n}$, $\mathbf{D_n}$, $\mathbf{E_6}$, $\mathbf{E_7}$,
$\mathbf{E_8}$, $\mathbf{F_4}$ and $\mathbf{G_2}$.

\section{The general construction}
 \label{sec:genconst}
In this section we describe a construction of principally
polarized abelian varieties admitting a given integral
representation of a finite group $G$.

\begin{remark}
Of course, every finite group has integral representations. Indeed, a
complex (irreducible) representation of degree $m$ defined over
the field $L$ gives rise to a rational (irreducible)
representation of degree $m \, [L : \mathbb{Q}]$, but any
(irreducible) rational representation is $\mathbb{Q}-$equivalent
to an integral representation.

Every finite group also has \textit{faithful} integral
representations. As an example, we have its regular representation.
\end{remark}

Recall that (see \cite{la-bi}) for any polarized abelian variety
$(X = V/\Lambda , H)$, a direct sum decomposition
$$
\Lambda = L_1 \oplus L_2
$$
is called a \textit{decomposition for $H$} if $L_1$ and $L_2$ are
isotropic with respect to $\Im H$. Equivalently, if there exists a
symplectic basis $\{\lambda_1 , \ldots, \lambda_n ,$ $\mu_1 ,
\ldots , \mu_n\}$ for $\Lambda$ such that $L_1 = \langle \lambda_1
, \ldots, \lambda_n \rangle$ and $L_2 = \langle \mu_1 , \ldots ,
\mu_n \rangle$.

We now define a similar notion for p.p.a.v.'s with $G-$action.

\begin{defi}
Let $G$ denote a finite group and consider a faithful integral
representation $\rho$ of $G$.

If $(X = V/\Lambda , H)$ is a p.p.a.v. with $G-$action, we will
say that a decomposition for $H$
$$
\Lambda = L_1 \oplus L_2
$$
is a $\rho-$\textit{decomposition} if both $L_1$ and $L_2$ are
$G-$invariant and $L_1$ is $\mathbb{Z}G-$equivalent to $\rho$.

In this case we will say that  $(X,H)$ is
$\rho-$\textit{decomposable}.
\end{defi}

\begin{remark}
If $(X  = V/\Lambda ,H)$ is $\rho-$decomposable, then there exists
a symplectic canonical basis for $\Lambda$ with respect to which
the $G-$action is represented by the symplectic matrices
\begin{equation} \label{eq:symp}
\left(%
\begin{array}{cc}
  \rho(g) & 0 \\
  0 & \rho^{-t}(g) \\
\end{array}%
\right) \ , \ g \in G
\end{equation}
where $\rho^{-t}$ is the contragredient representation of $\rho$.
The converse is also clearly true.

Also, the $G-$action on the tangent space to $X$ at the origin  is
(equivalent to) $\rho$.
\end{remark}

We will now give a construction for $\rho-$decomposable
p.p.a.v.'s.

\medskip

If $L$ is a $\mathbb{Z}G-$module corresponding to an integral
representation $\rho$, we let $U$ denote the real linear span of
$L$. For a $\rho(G)-$invariant real inner product $B$ on $U$ that
is integral on $L$, we will denote by $L^{*}$ the dual of $L$ with
respect to $B$.

\begin{proposition} \label{prop:exist}
Let $G$ denote a finite group and consider a faithful integral
representation $\rho$ of $G$ over a $\mathbb{Z}G$ module $L$.

Denote by $U$ (respectively $V$) the real (respectively complex)
linear span of $L$.

Then, for any $\rho(G)-$invariant real inner product $B$ on $U$
that is integral on $L$, and for any complex number $\tau$ with
positive imaginary part, the p.p.a.v. given by
$$
X = X(B, \tau) = V/(L + \tau L^*) \ , \ \Im H (u+\tau u' , v +
\tau v') = B(u,v') - B(v,u')
$$
is $\rho-$decomposable.

Furthermore, each such $X$ is isomorphic (as a torus) to a product
of elliptic curves.
\end{proposition}

\begin{proof}
Only the last assertion needs to be proven: since $L \subseteq
L^{*}$ are integral lattices of the same dimension $n$, then (see
\cite[p.\,221]{cass}) there exist a basis $\{ \ell_1 , \ldots ,
\ell_n\}$ of $L$ and positive integers $\{d_j\}_{j=1}^n$ such that
$d_{j+1}$ divides $d_j$ for each $j \in \{ 1 , \ldots , n-1\}$ and
such that $\{ \frac{1}{d_j} \, \ell_j\}_{j=1}^n$ is a basis for
$L^{*}$.

It is then clear that the torus $X = V/(L + \tau L^*)$ is
isomorphic to the product
$$
E_{\frac{\tau}{d_1}} \times \ldots \times E_{\frac{\tau}{d_n}}
$$
where $E_{\mu}$ denotes the elliptic curve with modulus $\mu$.
\end{proof}

When $\rho$ is an absolutely irreducible representation, we can
also prove the converse; that is, we now give a characterization
of $\rho-$de\-compo\-sable p.p.a.v.'s whenever $\rho$ is
absolutely irreducible.

\begin{theorem}
\label{theo:gentheo}
If in Proposition \ref{prop:exist} the
representation $\rho$ is absolutely irreducible, then every
$\rho-$decomposable p.p.a.v. is of the form described there.
\end{theorem}

\begin{proof}
Let $B$ denote a fixed real inner product on the real span of $L$,
with $B$ integral on $L$. For any $\tau$ in $\mathbb{H}_1$  we can
apply the construction of Proposition \ref{prop:exist} to find a
$\rho-$decomposable p.p.a.v. corresponding to $X(B, \tau)$.

Since $\rho$ is absolutely irreducible, any $\rho(G)-$invariant
inner product on the real linear span of $L_1$ is a positive
constant multiple of $B$ \cite[p.\,66]{b}; therefore, any positive
multiple $a \, B$ (integral on $L$) together with any $\tau$ in
$\mathbb{H}_1$ give rise (via the construction described in
Proposition \ref{prop:exist}) to the $\rho-$decomposable p.p.a.v.
corresponding to $X(a\, B, \tau)$, which is naturally isomorphic
to $X(B ,\frac{1}{a} \, \tau)$.

Given a $\rho-$decomposable p.p.a.v. $(X = V/\Lambda, H)$, without
loss of generality we assume that there exists a symplectic basis
$\{\lambda_1 , \ldots, \lambda_n ,$ $\mu_1 , \ldots , \mu_n\}$ for
$\Lambda$ such that $L_2 = \langle \mu_1 , \ldots , \mu_n \rangle$
and
$$
L = L_1 = \langle \lambda_1 , \ldots, \lambda_n \rangle .
$$

Then the Riemann matrix $Z = W+iY$ for $X$ with respect to this
basis is fixed under the symplectic action \eqref{eq:symp}, and
therefore $Y^{-1}$, the inverse of $Y$, yields a
$\rho(G)-$invariant inner product on the real linear span of $L$,
which is then a positive constant multiple of $B$.

Furthermore, $W \, Y^{-1}$ is a real endomorphism of the real
linear span of $L$ commuting with $\rho$, and therefore a multiple
of the identity.

We have thus proven that $Z$ is a complex multiple of the inverse
of the matrix of $B$ with respect to the basis $\lambda$, where
the scalar lies in $\mathbb{H}_1$.
\end{proof}

The next result follows immediately from the proof of the theorem.

\begin{corollary} \label{cor:riemannmat}
Let $G$ denote a finite group and consider a faithful integral
representation absolutely irreducible $\rho$ of $G$ on a
$\mathbb{Z}G-$module $L$.

Then the family of $\rho-$decomposable p.p.a.v.'s is parameterized
by the family of Riemann matrices
\begin{equation} \label{eq:Riemannmat}
Z_{\tau} = \tau Z_0 \ , \ \tau \in \mathbb{H}_1,
\end{equation}
where $Z_0$ is the inverse matrix of a real inner product on the
real span of $L$ which is integral on $L$, with respect to an
adequate basis for $L$.
\end{corollary}

\begin{remark}
Note that the Riemann matrices described in \eqref{eq:Riemannmat}
correspond to the fixed points in Siegel space under the action
given in \eqref{eq:symp}.
\end{remark}

Combining Proposition \ref{prop:exist} with Theorem
\ref{theo:gentheo} we obtain the following result.

\begin{corollary} \label{cor:absell}
Each $\rho-$decomposable p.p.a.v. for an absolutely irreducible
representation $\rho$ of $G$ is isomorphic (as a torus) to a
product of elliptic curves.
\end{corollary}

\begin{example}
The following is an example of an integral representation whose
associated family of p.p.a.v.'s is \textit{not} one-dimensional,
and where there are members of the family that are not invariant
under the natural action of any Weyl group in this dimension.

Consider the rational irreducible non-trivial representation of $G
= \mathbb{Z}/5\mathbb{Z} = \langle a : a^5\rangle$ given by
$$
 \rho (a) =
 \left(\begin{array}{rrrr}
           0 & 0 & 0 & -1 \\
           1 & 0 & 0 & -1 \\
           0 & 1 & 0 & -1 \\
           0 & 0 & 1 & -1
       \end{array}\right)
$$

Then the corresponding Riemann matrices are given by

\begin{equation}
\label{eq:z5}
 Z = \tau_1 \left(\begin{array}{rrrr}
4 & 3 & 2 & 1 \\
3 & 6 & 4 & 2 \\
2 & 4 & 6 & 3 \\
1 & 2 & 3 & 4
\end{array}\right)
+ \tau_2 \left(\begin{array}{rrrr}
2 & 2 & 1 & 0 \\
2 & 4 & 3 & 1 \\
1 & 3 & 4 & 2 \\
0 & 1 & 2 & 2
\end{array}\right)
\end{equation}

\noindent where $\tau_1$, $\tau_2$ are complex numbers such that
the imaginary part of $Z$ is positive definite.

Note that this representation is the restriction to the subgroup
$G$ of the natural representation of the Weyl group corresponding
to the root system $\mathbf{A_4}$, and that this is the only Weyl
group for an irreducible system of dimension $4$ having a subgroup
isomorphic to $G$. As we will see later on, the family of
p.p.a.v.'s associated to the natural representation of the Weyl
group just mentioned is the subfamily of (\ref{eq:z5}) obtained by
setting $\tau_2 = 0$.
\end{example}

\begin{example}
We give an example of an absolutely irreducible integral
representation acting on a one parameter family of p.p.a.v.'s of
dimension six (which contains a Jacobian, that of Wiman's curve,
see \cite{gonzr}), with the given action on the tangent space at
the origin of each member of the family, and such that there is an
invariant sublattice of rank six with \textit{no} invariant
complementary sublattice; that is, these p.p.a.v.'s are
\textit{not} $\rho-$decomposable.

The p.p.a.v.'s are defined by their Riemann matrices
$$
Z_{\tau} = \left[ {\begin{array}{cccccc} 3\,{\tau} &  - {\tau} &
{\tau} & - {\tau } & {\tau} & { -\frac{1}{2}}
\\ [2ex]
 - {\tau} & 3\,{\tau} &  - {\tau} &  - {\tau} & {\frac {1}{2}}  & {\tau} \\
[2ex] {\tau} &  - {\tau} & 3\,{\tau} & { -\frac{1}{2}}  & - {\tau}
& {\tau} \\ [2ex]
 - {\tau} &  - {\tau} & { -\frac{1}{2}}  & 3\,{\tau} &  - {\tau} &  - {\tau} \\
[2ex] {\tau} & {\frac {1}{2}}  &  - {\tau} & - {\tau} & 3\,{\tau}
&  - {\tau} \\
[2ex] {-\frac{1}{2}}  & {\tau} & {\tau} & - {\tau} &  - {\tau} &
3\,{\tau}
\end{array}}
 \right]
$$
for $\tau \in \mathbb{H}_1$.

The group is $G =S_5$, the symmetric group on five letters, and
the corresponding integral representation $\rho$ is given as
follows. Note that $\rho$ is the unique (complex) irreducible
representation of degree six of $S_5$.

$$
\rho (1,2,3,4,5) = \left[ {\begin{array}{rrrrrr}
1 & 1 & 0 & 1 & 0 & 0  \\
-1 & 0 & 1 & 0 & 1 & 0  \\
1 & 0 & 0 & 0 & 0 & 0  \\
0 & -1 & -1 & 0 & 0 & 1  \\
0 & 1 & 0 & 0 & 0 & 0  \\
0 & 0 & 1 & 0 & 0 & 0
\end{array}}
 \right]
 $$
$$
\rho (1,2)=  \left[ {\begin{array}{rrrrrr}
1 & 0 & 0 & 0 & 0 & 0  \\
0 & 1 & 0 & 0 & 0 & 0  \\
0 & 0 & 1 & 0 & 0 & 0  \\
-1 & -1 & 0 & -1 & 0 & 0  \\
1 & 0 & -1 & 0 & -1 & 0  \\
0 & 1 & 1 & 0 & 0 & -1
\end{array}}
 \right]
$$

The associated symplectic action is given by
$$
\widetilde{\rho} (1,2,3,4,5) =
   \left[ {\begin{array}{cc} \rho (1,2,3,4,5) & 0 \\
                             0 & \rho^{-t} (1,2,3,4,5)
            \end{array}}
   \right]
$$
$$
    \widetilde{\rho} (1,2)=
       \left[ {\begin{array}{cc}  \rho (1,2) & L(1,2) \\
                             0 & \rho^{-t} (1,2)
            \end{array}}
   \right]
$$
where
$$
L(1,2) =
 \left[ {\begin{array}{rrrrrr}
 0 & 0 & 0 & 0 & 0 & 1 \\
 0 & 0 & 0 & 0 & -1 & 0 \\
 0 & 0 & 0 & 1 & 0 & 0 \\
 0 & 0 & -1 & 0 & 1 & -1 \\
 0 & 1 & 0 & -1 & 0 & 1 \\
 -1 & 0 & 0 & 1 & -1 & 0
 \end{array}}
   \right]
$$

The Riemann matrices $Z_{\tau}$ are the only Riemann matrices
invariant under the given action, and a calculation shows that the
given symplectic representation is not even integrally equivalent
to any integral representation of the following form
$$
\left(%
\begin{array}{cc}
  U(g) & 0 \\
  0 & U(g)^{-t} \\
\end{array}%
\right) \ .
$$
\end{example}

\section{The construction for Weyl groups}
\label{sec:constweyl}

For each Weyl group $\mathcal{W}$ there is a natural integral
faithful representation $\rho$, generated by the reflections on
the roots. Furthermore, this representation is absolutely
irreducible if the associated root system is irreducible and
reduced.

Therefore, in this case we may apply Theorem~\ref{theo:gentheo} to
find the family of $\rho-$decomposable p.p.a.v.'s, which we also
know is parameterized by $\mathbb{H}_1$ .

We summarize this information in the next proposition, where we
include the Riemann matrices for each family for the sake of
completeness.

\begin{proposition}
 \label{prop:theoweyl}
Let $R$ be an irreducible, reduced root system of dimension $n$,
and let $\mathcal{W}(R)$ be its Weyl group.

We denote the corresponding canonical basis for the root system
$(V,R)$ by $\mathcal{B}$, and consider the integral representation
$\rho$ of dimension $n$ of $\mathcal{W}(R)$ generated by the
reflections associated to the elements of $\mathcal{B}$

Then the family of $\rho-$decomposable p.p.a.v.'s is
one-dimensional, and parameterized by the family of Riemann
matrices $\, \mathbb{H}(R) \,$ in $\mathbb{H}_n$ associated to
$\rho$ by Corollary~\ref{cor:riemannmat}.

The corresponding Riemann matrices are given by

$$
\mathbb{H}(R) = \{ Z_{\tau} = \tau Z_0 : \tau \in \mathbb{H}_1 \}
 $$

\smallskip

\noindent where $Z_0$ is given in Table \ref{table:z0}.
\end{proposition}

\vspace{5mm}

The proof of Proposition~\ref{prop:theoweyl} follows Table
\ref{table:z0}. Note that since $Z_0$ is symmetric, we give only
the upper half of the matrix.

\newpage

 {\scriptsize{
\begin{center}
\begin{longtable}{|c|c|}
\caption[]{The Riemann matrices
\label{table:z0}}\\
\hline
 \endfirsthead
 \caption[]{(continued)}\\ \hline
 \endhead
             &                                  \\
 Root system & $Z_0$                  \\*
             &                            \\* \hline
             &                            \\
 $\mathbf{A_n} $      & $\dfrac{1}{n+1} \left(\begin{array}{ccccccc}
n & n-1    & n-2      & \ldots  &        &  2  & 1 \\
  & 2(n-1) &          & \ldots  &        &  4  & 2 \\
  &        & \ddots   &         &        &     & \\
  &        &          & j(n+1-j)& \ldots & 2j  &  j \\
  &        &          &         & \ddots &     &  \\
  &        &          &         &        & 2(n-1) & n-1 \\
  &        &          &         &        &     & n
\end{array}
\right)$      \\*
             &                                       \\ \hline
             &                                       \\*
 $\mathbf{B_n} $      & $\left(\begin{array}{ccccc}
1  & 1 & \ldots  & 1  & 1 \\
 & 2 & \ldots  & 2  & 2 \\
   &   & \ddots  &    & \\
  &  &   & n-1  & n-1 \\
  &  &   &   & n
\end{array}
\right)$         \\*
             &                                     \\ \hline
             &                                     \\*
 $\mathbf{C_n}$       & $ \left(\begin{array}{ccccc}
1  & 1 & \ldots  & 1 & \frac{1}{2} \\*
 & 2 & \ldots  & 2   & 1 \\
   &   & \ddots  &    & \\
  &  &   & n-1  & \frac{n-1}{2} \\
  &  &   &   & \frac{n}{4}
\end{array}
\right) $       \\*
             &                                       \\ \hline
             &                                       \\*
 $\mathbf{D_n}$       &     $\left(\begin{array}{cccccccc}
1  & 1 & 1       & 1 & 1      & \ldots &\frac{1}{2} &
\frac{1}{2}\\*
   & 2 & 2       & 2 & 2      & \ldots & 1          & 1 \\*
   &   & \ddots  &   &        &        &            & \\*
   &   &         & j & j      & \ldots & \frac{j}{2} &
   \frac{j}{2}\\*
   &   &         &   & \ddots &        &             &
   \\*
   &   &         &   &        &  n-2   & \frac{n-2}{2}& \frac{n-2}{2}
   \\*
   &   &         &   &        &        & \frac{n}{4}& \frac{n-2}{4}
   \\*
   &   &         &   &        &        &              & \frac{n}{4}
\end{array}
\right) $        \\*
            &             \\ \hline
           &             \\*
 $\mathbf{E_6}$       & $\left( \begin{array}{rrrrrr}
\frac{4}{3} & 1 & \frac{5}{3}  & 2 & \frac{4}{3} & \frac{2}{3} \\
            & 2 & 2            & 3 & 2           & 1           \\
            &   & \frac{10}{3} & 4 & \frac{8}{3} & \frac{4}{3} \\
            &   &              & 6 & 4           & 2           \\
            &   &              &   & \frac{10}{3}& \frac{5}{3}\\
            &   &              &   &             & \frac{4}{3}
       \end{array}
       \right)$        \\
             &                                      \\ \hline
             &                                      \\*
 $\mathbf{E_7}$      & $\left( \begin{array}{rrrrrrr}
                    2  & 0  & -1 & 0  & 0  & 0  & 0  \\*
                       & 2  & 0  & -1 & 0  & 0  & 0 \\
                       &    & 2  & -1 & 0  & 0  & 0  \\
                       &    &    & 2  & -1 & 0  & 0 \\
                       &    &    &    & 2  & -1 & 0 \\
                       &    &    &    &    & 2  & -1\\
                       &    &    &    &    &    & 2
                    \end{array}
       \right)$      \\*
             &                                      \\ \hline
             &                                     \\*
 $\mathbf{E_8}$       & $\left( \begin{array}{rrrrrrrr}
                    4  & 5  & 7  & 10 & 8  & 6  & 4  & 2 \\
                       & 8  & 10 & 15 & 12 & 9  & 6  & 3\\
                       &    & 14 & 20 & 16 & 12 & 8  & 4 \\
                       &    &    & 30 & 24 & 18 & 12 & 6\\
                       &    &    &    & 20 & 15 & 10 & 5\\
                       &    &    &    &    & 22 & 8  & 4\\
                       &    &    &    &    &    & 6  & 3\\
                       &    &    &    &    &    &    & 2
                    \end{array}
       \right)$           \\*
             &                                    \\ \hline
             &                                   \\*
 $\mathbf{F_4}$       & $\left(
       \begin{array}{rrrr}
    1 & \frac{3}{2} & 2 & 1 \\
      & 3 & 4 & 2 \\
      &   & 6 & 3 \\
      &   &   & 2
       \end{array} \right)$      \\*
             &                                       \\ \hline
             &                                      \\*
 $\mathbf{G_2}$       & $\dfrac{1}{3} \left( \begin{array}{cc}
            6 & 3 \\
             & 2 \end{array} \right)$         \\
             &                                       \\ \hline
  \end{longtable}
\end{center}
}}
\normalsize

\begin{proof}[\textit{Proof of Proposition~\ref{prop:theoweyl}.}]

We follow the notation of \cite{b}; especially see Planches I
through IX. All matrices will be given with respect to the basis
$\mathcal{B} = \{ \alpha_1 , \ldots , \alpha_n \}$.

Let $C$ denote the Cartan matrix of $R$ given by
$$ C = \left( c_{ij} \right) =
\left( 2 \frac{( \alpha_i , \alpha_j)}{(\alpha_j,\alpha_j)}
\right) \ ,
$$
\noindent and let $D$ denote the $n \times n$ diagonal matrix
given by $d_{ii} = \dfrac{(\alpha_i , \alpha_i)}{2} $.

\medskip

Then an integral matrix $S$ for the $\mathcal{W}(R) - $invariant
inner product is given by $S = CD$, which allows us to compute the
values for $Z_0 = S^{-1}$ given in the table.
\end{proof}

\begin{remark}
Let $R$ be an irreducible, reduced root system of dimension $n$.
It is known (\cite{b}) that the group of automorphisms of $R$,
denoted by $\aut(R)$, is a semidirect product of the corresponding
Weyl group $\mathcal{W}(R)$ and a finite group $\Gamma$, where
$\Gamma$ is the trivial group except in the following cases.

\begin{description}
  \item[$\mathbf{R = A_n} \,$] $\Gamma = \mathbb{Z}/2\mathbb{Z}$

\medskip

\item[$\mathbf{R = D_4} \,$] $\Gamma = S_3$

\medskip

\item[$\mathbf{R = D_n, n > 4} \,$] $\Gamma =
\mathbb{Z}/2\mathbb{Z}$

\medskip

\item[$\mathbf{R = E_6}\,$] $\Gamma = \mathbb{Z}/2\mathbb{Z}$
\end{description}
\end{remark}

A case by case computation allows us to prove the following
result.

\begin{corollary}
\label{cor:autr}
Let $R$ be an irreducible, reduced root system of
dimension $n$.

Denote by $\aut(R)$ the group of automorphisms of $R$ and by
$\widetilde{\rho}$ the natural integral representation of
$\aut(R)$ extending the representation $\rho$ of the Weyl group
$\mathcal{W}(R)$.

Then the family of $\widetilde{\rho}-$decomposable p.p.a.v.'s
coincides with the family of $\rho-$decomposable p.p.a.v.'s.

More explicitly, each $Z_{\tau} \in \mathbb{H}(R)$ is fixed under
the symplectic action \eqref{eq:symp} induced by
$\widetilde{\rho}$.
\end{corollary}

\medskip

\begin{defi} \label{defi:vapp}
Let $R$ be an irreducible, reduced root system of dimension $n$,
and let $\mathcal{W}(R)$ be its Weyl group.

Consider the family  $\, \mathbb{H}(R) \,$ of Riemann matrices in
$\mathbb{H}_n$ associated to $\rho(\mathcal{W}(R))$ by Proposition
\ref{prop:exist}.

For each $Z_{\tau}$ in $\, \mathbb{H}(R) \,$ we will denote by
$A_{\tau} = (X_{\tau}, H_{\tau})$ the corresponding
$\rho-$decomposable p.p.a.v., and by $\mathcal{A}(R)$ the
corresponding family in $\mathcal{A}_n$.

Note that each $A_{\tau}$ is $\widetilde{\rho}-$decomposable, for
$\widetilde{\rho}$ as in Corollary \ref{cor:autr}.
\end{defi}

\begin{proposition}
\label{prop:isoaut}
With the notation of Proposition
\ref{prop:theoweyl}, the following results hold for each $\tau$ in
$\mathbb{H}_1$, where the corresponding isomorphisms are as
p.p.a.v.'s.

\medskip

\begin{enumerate}
\item $A_{\tau} \in \mathbb{H}(\mathbf{D_n})$ is isomorphic to $
A_{\tau} \in \mathbb{H}(\mathbf{C_n})$ for each $n$,

\vskip12pt

\item $A_{\tau} \in \mathbb{H}(\mathbf{F_4})$ is isomorphic to
$A_{\tau} \in \mathbb{H}(\mathbf{D_4})$, and

\vskip12pt

\item $A_{\tau} \in  \mathbb{H}(\mathbf{G_2})$ is  isomorphic to
$A_{\tau} \in \mathbb{H}(\mathbf{A_2})$.
\end{enumerate}
\end{proposition}

\begin{proof}
It is easy to verify that, in each of the above cases, the
integral representations $\widetilde{\rho}$ considered in
Corollary \ref{cor:autr} are isomorphic over $\mathbb{Z}$ for the
corresponding groups $\aut(R)$.
\end{proof}

\begin{remark}
Alternatively, the results of Proposition \ref{prop:isoaut} may be
proved by explicit computations: to show each member of a family
$F_1 =\mathbb{H}(\mathbf{R_1})$ is isomorphic to the corresponding
element of the family $F_2 =\mathbb{H}(\mathbf{R_2})$, we find an
invertible complex $n \times n$ matrix $A$ and a symplectic $2n
\times 2n$ matrix $M$ such that the following equation holds for
each $Z_{\tau} = \tau Z_1$ in $F_1$ and each $\widetilde{Z_{\tau}}
= \tau Z_2$ in $F_2$.

\begin{equation}
\label{eq:iso}
 A \left( I_n \ \ Z_{\tau} \right) = \left( I_n \ \
\widetilde{Z_{\tau}} \right) M \ .
\end{equation}

In fact, in each case we can choose $A$ to be unimodular
(integral) and $M$ to be of the form
 $M = \left(
      \begin{smallmatrix}
      A & 0 \\
      0 & \mbox{}^t A^{-1}
      \end{smallmatrix}
      \right) \ .
 $

For $F_1 = \mathbb{H}(\mathbf{D_n})$ and $F_2 =
\mathbb{H}(\mathbf{C_n})$ we let $A = \left(
\begin{array}{cccr}
        \multicolumn{2}{c}{I_{n-2}}  & \multicolumn{2}{c}{0} \\
        0 & 1 & 1 & -1 \\
        0 & 0 & 1 & 0
\end{array}
\right)$.

For $F_1 = \mathbb{H}(\mathbf{D_4})$ and $F_2 =
\mathbb{H}(\mathbf{F_4})$ we let $A = \left(
\begin{array}{cccc}
        0 & 0 & 0 & 1 \\
        1 & 0 & 0 & 1 \\
        1 & 0 & 1 & 1 \\
        0 & 1 & 0 & 0
\end{array}
\right)$.

For $F_1 = \mathbb{H}(\mathbf{G_2})$ and $F_2 = \mathbb{H}
(\mathbf{A_2})$ we let $A = \left( \begin{array}{cr}
        1 & -2 \\
        0 & 1
\end{array}
\right)$.
\end{remark}

\begin{remark}
It will be a consequence of Proposition \ref{prop:dec}, Corollary
\ref{cor:bn} and Theorem \ref{theo:pavirred} that there no other
isomorphisms between different families of type
$\mathbb{H}(\mathbf{R})$.
\end{remark}

\begin{remark}
\label{rem:an} The family $\mathbb{H}(\mathbf{A_n})$ was discussed
in a different context in \cite{gr}. The corresponding Riemann
matrices mentioned there are of the following form

$$ \widetilde{Z}_{\tau} =
\tau \left(
\begin{array}{rrrr}
n & - 1 &  \ldots & -1\\ -1 & n & & \vdots\\ & & \ddots &  -1\\ -1
&  \cdots & -1 & n
\end{array} \right) \ .
$$

This family is of course isomorphic to our family; as above, it is
enough to consider the unimodular $n \times n$ matrix

$$ A = \left(
\begin{array}{rrrrr}
1 & 0 & 0 & 0 & \cdots\\
-1 & 1 & 0 & 0 & \cdots \\
0 & -1 & 1 & 0  & \cdots \\
\vdots &  & \ddots & \ddots & \\
0 &   \ldots & & -1 & 1
\end{array} \right)
$$

\noindent and verify that Equation (\ref{eq:iso}) is satisfied
with $M = \left(
      \begin{smallmatrix}
      A & 0 \\
      0 & A^{-t}
      \end{smallmatrix}
      \right)
 $.
\end{remark}

\begin{remark}
 \label{rem:loo}
If $R$ is an irreducible, reduced root system we can also perform
the construction in Theorem \ref{theo:gentheo} for the natural
integral representation associated to the lattice generated by
$R^{\vee}$, the inverse root system.

A related construction is the following. Let $E$ be a generic
elliptic curve and $R$ an irreducible, reduced root system of rank
$n$, with $Q$ the lattice generated by $R$ or $R^{\vee}$. Then $X
\simeq Q\otimes E$ is an $n$-dimensional complex torus isomorphic
(as a complex torus) to the $n$-fold product of the elliptic curve
$E$ and the Weyl group $\mathcal{W}(R)$ acts naturally on $X$. The
polarization $P$ on $E$ induces an equivariant polarization $P_X$
on $X$ which is, in general, not principal.

Such abelian varieties have been discussed earlier in \cite{loo}
for the case $Q^{\vee}$ the lattice generated by $R^{\vee}$. One
can verify that the natural equivariant polarization induced on $X
\simeq Q^{\vee}\otimes E$ is of respective degree $n+1$, $4$, $1$,
$4$, $3$, $2$, $1$, $4$ and $3$, corresponding to the system
$\mathbf{A_n}$, $\mathbf{B_n}$, $\mathbf{C_n}$, $\mathbf{D_n}$,
$\mathbf{E_6}$, $\mathbf{E_8}$, $\mathbf{F_4}$ and $\mathbf{G_2}$.

In Section \ref{sec:modular} we will make explicit the relation
between these polarized abelian varieties and our corresponding
principally polarized abelian varieties.
\end{remark}

\section{Irreducibility of the abelian varieties associated to Weyl groups}
\label{sec:irred} In this section we answer the question of
irreducibility of the p.p.a.v.'s admitting the natural action of
the Weyl groups constructed in Proposition \ref{prop:theoweyl}.

We already know, from Corollary \ref{cor:absell}, that  when
considered as complex tori they are all isomorphic to the product
of elliptic curves. Our next result gives the explicit
decomposition in each case.

\begin{proposition}
\label{prop:dec} Let $R$ be an irreducible, reduced root system
and consider the corresponding family of p.p.a.v.'s
$\mathcal{A}(R)$ as per Definition \ref{defi:vapp}.

Then each $A_{\tau}$ in $\mathcal{A}(R)$ is isomorphic (as a
complex torus) to a product of elliptic curves, as follows.

\newpage

\small{
\begin{center}
\begin{longtable}{|c|c|c|}
\caption[]{Decomposition as complex tori
\label{table:tori}}\\
\hline
 \endfirsthead
 \caption[]{(continued)}\\ \hline
 \endhead
             &                               &   \\*
 Root system & Decomposition for $A_{\tau}$                & Case \\*
             &                                &      \\ \hline
             &                                &     \\*
 $\mathbf{A_n} $      & $E_{\tau}^{n-1}  \times E_{\tau/(n+1)}$  & any $n$ \\*
             &                                 &      \\ \hline
             &                                 &      \\*
 $\mathbf{B_n} $      & $E_{\tau}^{n}$     & any $n$    \\*
             &                                 &      \\ \hline
             &                                 &      \\*
 $\mathbf{C_n}$, $\mathbf{D_n}$& $E_{\tau}^{n-1}  \times E_{\tau/4}$   &  $n$ odd \\*
             &                                             & \\*
             & $E_{\tau}^{n-2}  \times E_{\tau/2}^{\, 2}$  & $n$ even \\*
             &                                 &      \\ \hline
             &                                 &      \\*
 $\mathbf{E_6}$       & $E_{\tau}^{5} \times E_{\tau/3}$  & \\*
             &                                 &      \\ \hline
             &                                 &     \\*
 $\mathbf{E_7}$       & $E_{\tau}^{6}  \times E_{\tau/2}$  & \\*
             &                                 &      \\ \hline
             &                                 &     \\*
 $\mathbf{E_8}$       & $E_{\tau}^{8}$      &     \\*
             &                                 &      \\ \hline
             &                                 &     \\*
 $\mathbf{F_4}$       & $E_{\tau}^2 \times E_{\tau/2}^2$      & \\*
             &                                 &      \\ \hline
             &                                 &     \\*
 $\mathbf{G_2}$       & $E_{\tau} \times E_{\tau/3}$      & \\*
             &                                 &      \\ \hline
 \end{longtable}
\end{center}
}
\end{proposition}

\begin{proof}
A case by case calculation of the integers $d_j$ appearing in the
proof of Proposition \ref{prop:exist} gives the result.
\end{proof}

\begin{remark} \label{rem:bn}
Alternatively, given a Riemann matrix $Z_{\tau}$ in
$\mathbb{H}(R)$, we can find explicit $n \times n$ unimodular
matrices $F$ and $M$ and a diagonal  $n \times n$ matrix $d$ (with
positive entries in the diagonal) which satisfy the following
equation

\begin{equation}
\label{eq:dec} F \ ( I_n \ \  Z_{\tau} ) = (I_n \ \ \tau d)\left(
\begin{array}{cc}
F & 0 \\
0 & M
\end{array} \right)
\end{equation}

We illustrate with the following example, which will give extra
information.

\vskip12pt

Let $R =  \mathbf{B_n}$. In this case, let
 $$
 F = \left(
\begin{array}{rrrrc}
                    1 & 0 & \ldots & 0 & 0 \\
                   -1 & 1 & \ldots & 0 & 0 \\
                      &   & \ddots &   & \\
                    0 &   & \ldots & 1 & 0 \\
                    0 &   & \ldots & -1& 1 \\
                    \end{array}
       \right) \ , \ d = I_n \ , \ \text{ and } \ M = FS^{-1} \ ,
       $$
and note that equation (\ref{eq:dec}) is satisfied.
\end{remark}

\medskip

This calculation provides the first result about the
irreducibility of these p.p.a.v.'s, as follows.

\begin{corollary}
\label{cor:bn}
For every natural number $n$ and every $\tau$ in
$\mathbb{H}_1$, the p.p.a.v. $A_{\tau}$  in
$\mathcal{A}(\mathbf{B_n})$ is isomorphic as a p.p.a.v. to the
product of $n$-times the elliptic curve $E_{\tau}$.
\end{corollary}

\begin{proof}
Observe that $M = F^{-t}$ for the explicit matrices given above.

Therefore the matrix
 $\left(
\begin{array}{cc}
F & 0 \\ 0 & M
\end{array} \right)$
appearing in the equivalence of Remark \ref{rem:bn} is symplectic,
which finishes the proof.
\end{proof}

\bigskip

About the case $R =  \mathbf{A_n}$, the following result may be
found in \cite{gr} (see also Remark \ref{rem:an}).

\begin{proposition}
Consider the root system of type $\mathbf{A}_{n}$ and denote by
$\mathcal{A}(\mathbf{A_{n}})$ the corresponding family of
p.p.a.v.'s of dimension $n$ given in Definition \ref{defi:vapp}.

Then each $A_{\tau}$ in $\mathcal{A}(\mathbf{A_n})$ is irreducible
as a p.p.a.v., except when $n=2$ and $\tau$ is equivalent to $
\tau_0 = \dfrac{-3 + i\sqrt{3}}{6}$ mod $\Gamma_0 (3)$, in which
case $A_{\tau}$ is isomorphic as a p.p.a.v. to $E_{\omega}\times
E_{\omega}$, where $E_{\omega}$ is the elliptic curve with
$j(\omega)= 0$.
\end{proposition}

\medskip

We will now answer the question of irreducibility as p.p.a.v.'s
for the other irreducible root systems, where by the isomorphisms
given in Proposition \ref{prop:isoaut} we only need to consider
the cases $R = \mathbf{C_n}$, $\mathbf{E_6}$, $\mathbf{E_7}$ and
$\mathbf{E_8}$. Note that the proof applies to the case
$\mathbf{A}_n$ as well.

\begin{theorem}
\label{theo:pavirred} For each of the root systems $R =
\mathbf{C_n}$ ($n \geq 3$), $\mathbf{E_6}$, $\mathbf{E_7}$ and
$\mathbf{E_8}$ consider the respective family of p.p.a.v.'s
$\mathcal{A}(R)$ given in Definition \ref{defi:vapp}.

Then each $(X_{\tau}, H_{\tau})$ in $\mathcal{A}(R)$  is
irreducible as a p.p.a.v., except for some exceptional values of
$\tau$ in the case $R = \mathbf{C_n}$.
\end{theorem}

\begin{proof}
Suppose that $(X_{\tau}, H_{\tau}) = (X,H)$ is isomorphic as a
p.p.a.v. to a product of p.p.a.v.'s as follows
$$  (X, H) \simeq   (X_1 , H_1) \times  (X_2 , H_2) \times \ldots \times (X_k ,
H_k) ,
$$
with $k > 1$.

Then the tangent space at the origin for $X$, $T(X)$, decomposes
accordingly
$$
T(X) \simeq  \oplus_{j=1}^k T(X_j) \ .
$$

Furthermore, the corresponding Weyl group preserves this
decomposition, and therefore the group is imprimitive.

But then the root system under consideration must be
$\mathbf{C_n}$ (see \cite{gm}); furthermore all $(X_j,H_j)$ are
isomorphic to each other, and $\dim X_j = 1$.

That is, $X \: \simeq \:  X_1^n $, where $X_1$ is an elliptic
curve. But by Proposition \ref{prop:dec} we also have that
$$ X  \: \simeq \:  E_{\tau}^{n-1}  \times E_{\tau/4}
\; \; \; \mbox{ for $\; n \; $ odd} $$
 \noindent and that
$$  X \: \simeq \: E_{\tau}^{n-2}  \times E_{\tau/2}
  \times E_{\tau/2} \; \; \; \text{ for $ \; n \; $ even}$$
\noindent as complex tori.

This finishes the proof, since then the only values of $\tau$ for
which $(X_{\tau}, H_{\tau})$ may be reducible are for the case of
the root system $\mathbf{C_n}$ and $E_{\tau} \simeq E_{\tau/4}$,
for $n$ odd, and $E_{\tau} \simeq E_{\tau/2}$, for $n$ even.
\end{proof}

\begin{remark}
The Jacobian  variety $J(C)$ of a complex curve  $C$ of genus $n$
is an irreducible principally polarized abelian variety of
dimension $n$ (irreducibility of the theta divisor). Nevertheless,
it is an interesting problem to find Jacobians isomorphic as tori
(and not only isogenous) to products of elliptic curves.

As all the p.p.a.v.'s in the families associated to Weyl groups
are \textit{isomorphic} to products of elliptic curves, a
necessary condition for a  family $\mathcal{A}(R)$ in the moduli
space $\mathcal{A}_n$ of p.p.a.v.'s of dimension $n$ to intersect
the Jacobian locus ${\mathcal J}_n$ is given by the Hurwitz and
Torelli theorems.

For example, these condition are satisfied for the root systems
$\mathbf{A_4}$ and $\mathbf{A_5}$, and in those cases the
corresponding intersection is characterized  in \cite{rg},
\cite{gr} and \cite{gr2}.
\end{remark}

\section{The abelian varieties families for Weyl groups are modular curves}
\label{sec:modular}

In this section we will show that the families $\mathcal{A}(R)$
constructed in Section \ref{sec:constweyl} may be described as
modular curves in the corresponding moduli space of principally
polarized abelian varieties.

\begin{remark}
Let $R$ be an irreducible root system and denote by $G$ the
natural symplectic representation of either its Weyl group
$\mathcal{W}(R)$ given by \eqref{eq:symp}.

Denote by $C(G)$ and $N(G)$ the respective centralizer and
normalizer of $G$ inside the corresponding symplectic group, and
by $\aut(G)$ the group of automorphisms of $G$ that may be
realized inside the symplectic group.

Then  we have the following exact sequence of groups
 $$ 1 \to C(G) \overset{i}{\to} N(G) \overset{\varphi}{\to} \aut(G) \to 1
 $$
where $i$ denotes the inclusion and $\varphi(M)$ is conjugation by
$M$ for each $M$ in $ N(G)$. Furthermore, the family
$\mathcal{A}(R)$ in $\mathcal{A}_n$ is given by
$\mathbb{H}(R)/C(G)$, so we need to describe $C(G)$, which we do
next.
\end{remark}

\begin{proposition}
 \label{prop:centr}
 Let $R$ be an irreducible root system of dimension $n$
and consider the symplectic representation of its Weyl group
$\mathcal{W}(R)$ associated to the natural basis for the roots, as
in \eqref{eq:symp}.

Then its centralizer in $\spp(2n, {\mathbb Z})$ is isomorphic to a
subgroup of  $\psl (2,{\mathbb Z}) = \Gamma$, as follows.

\small{\begin{center}
\begin{longtable}{|c|c|c|}
\caption[]{Centralizers for $\mathcal{W}(R)$
\label{table:cent}}\\
\hline
 \endfirsthead
 \caption[]{(continued)}\\ \hline
 \endhead
             &                                 &     \\*
 Root system & Centralizer        & Case \\*
             &                                 &       \\* \hline
             &                                 &      \\*
 $\mathbf{A_n}$       & $\Gamma^0(n+1)$ & any $n$    \\*
             &                                 &      \\* \hline
             &                                 &      \\*
 $\mathbf{B_n}$, $\mathbf{E_8}$& $\Gamma$     & any $n$    \\*
             &                                 &      \\ \hline
             &                                 &      \\*
 $\mathbf{C_n}$, $\mathbf{D_n}$& $\Gamma^0(4)$  &  $n$ odd    \\*
             & $\Gamma^0(2)$  & $n$ even \\*
             &                                 &      \\ \hline
             &                                 &      \\*
 $\mathbf{E_6}$, $\mathbf{G_2}$& $\Gamma^0(3)$   &      \\*
             &                                 &      \\ \hline
             &                                 &     \\*
 $\mathbf{E_7}$, $\mathbf{F_4}$& $\Gamma^0(2)$   &     \\*
             &                                 &      \\ \hline
 \end{longtable}
\end{center}
}
\end{proposition}

\begin{proof}
We observe that $\Gamma = \psl(2, \mathbb{Z})$ acts naturally on
$\mathbb{Z}^2 \otimes L^{*} = L^{*} \oplus L^{*}$, and that the
centralizer we are looking for is isomorphic to the stabilizer of
$L \oplus L^{*}$ under this action, which is then computed using
Table \ref{table:tori}.
\end{proof}

\begin{remark}
Alternatively, since the natural integral representation of each
Weyl group is absolutely irreducible, it follows that each matrix
in the centralizer has the following form

 $$ \left( \begin{array}{cc}
         a I_n & b Z_0  \\
         c Z_0^{-1} & d I_n
          \end{array}
    \right)
 $$
with $a, b, c, d \in {\mathbb Z}$ and $ad - bc = 1$, and $Z_0$ as
in Table \ref{prop:theoweyl}. Noting that $Z_0^{-1}$ is an
integral matrix, the specific extra conditions on $b$ in each case
follow from Table \ref{prop:theoweyl}.
\end{remark}

\medskip

We can now parameterize each family $\mathcal{A}(R)$ as follows.

\begin{theorem}
\label{theo:mod} Let $R$ be an irreducible root system of
dimension $n$.

Then  the set $\mathcal{A}(R)$ of isomorphism classes of
p.p.a.v.'s in $\mathcal{A}_n$ which are $\rho-$decomposable for
the natural action of the corresponding Weyl group
$\mathcal{W}(R)$, is parameterized by a modular curve.

More precisely, we have the following result.

\small{
\begin{center}
\begin{longtable}{|c|c|c|}
\caption[]{Parameterizing the families $\mathcal{A}(R)$
\label{table:modular}}\\
\hline
 \endfirsthead
 \caption[]{(continued)}\\ \hline
 \endhead
             &                                 &     \\*
 Root system $R$&     $\mathcal{A}(R)$             & Case \\*
             &                                 &       \\* \hline
             &                                 &      \\*
 $\mathbf{A_n} $      & ${\mathbb H}_1/\Gamma^0(n+1)$ & any $n$    \\
             &                                 &      \\ \hline
             &                                 &      \\*
 $\mathbf{B_n}$, $\mathbf{E_8}$& ${\mathbb H}_1/\Gamma$     & any $n$ \\*
             &                                 &      \\ \hline
             &                                 &      \\*
 $\mathbf{C_n}$, $\mathbf{D_n}$&  ${\mathbb H}_1/\Gamma^0(4)$  &  $n$ odd \\*
             &                                 & \\*
             &  ${\mathbb H}_1/\Gamma^0(2)$  & $n$ even \\*
             &                                 &      \\ \hline
             &                                 &      \\*
 $\mathbf{E_6}$, $\mathbf{G_2}$& ${\mathbb H}_1/\Gamma^0(3)$   & \\*
             &                                 &      \\ \hline
             &                                 &     \\*
 $\mathbf{E_7}$, $\mathbf{F_4}$& ${\mathbb H}_1/\Gamma^0(2)$   & \\*
             &                                 &      \\ \hline
\end{longtable}
\end{center}
}
\end{theorem}

\begin{proof}
Consider an irreducible root system $R$ of dimension $n$ and
denote by $\mathbb{H}(R)$ the family of Riemann matrices in
$\mathbb{H}_n$ constructed in Proposition \ref{prop:theoweyl}.

It follows that the natural morphism
 $$ \mathbb{H}(R) \to \mathbb{H}_1
 $$
which sends $Z_{\tau}$ to $\tau$ is equivariant under the action
of the centralizer of $\aut(R)$ in $\spp(2n, {\mathbb Z})$.

Since this centralizer was described in Proposition
\ref{prop:centr} as a subgroup $\Gamma_1$ of $\psl(2,{\mathbb
Z})$, there is an induced injective morphism
 $$ \mathbb{H}(R)/\text{centralizer} \to
 \mathbb{H}_1/\Gamma_1 \, .
 $$

Conversely, let $(E,U)$ be a point of the modular curve
$\mathbb{H}_1/\Gamma_1$. Then $E$ is an elliptic curve and $U$ is
a cyclic subgroup of the appropriate order in $E$.

The construction of a principally polarized abelian variety in our
family starting from these data is obtained following the idea
discussed in \cite{gr} for the case $\mathbf{A_n}$ as follows.
Start by constructing an abelian variety $X$ of dimension $n$,
naturally associated to the root system and the elliptic curve,
and with the natural action of the Weyl group on $X$. There is a
natural equivariant polarization $P_X$ on $X$ such that we can
embed the subgroup $U$ inside the kernel of $P_X$, obtaining a
cyclic group which is both invariant under the group action and
totally isotropic for the Weil form. It follows that there is a
principal polarization $Q$ on $X/U$ whose inverse image is $P_X$.

Since the construction is equivariant for the group action, the
principally polarized abelian variety $(X/U,Q)$ is the required
one.

The proof is finished by explicit calculations in a case by case
analysis, by giving an embedding of each cyclic group $U$ into the
kernel of $P_X$ with the required properties.
\end{proof}

\begin{remark}
\label{rem:loo2} We can now complete the description of the
relationship between our varieties and the ones constructed in
\cite{loo}.

With the notation given in Remark \ref{rem:loo}, the abelian
varieties $X = Q^{\vee}\otimes E$ of \cite{loo} coincide with the
abelian varieties $X$ constructed in the proof of Theorem
\ref{theo:mod} (for the respective inverse root system). In other
words, our p.p.a.v.'s are obtained as quotients by appropriate
cyclic subgroups of the abelian varieties given in \cite{loo}.
\end{remark}

\begin{ack}
The authors would like to thank the referee for making numerous
suggestions that improved the presentation and for pointing out
misprints.
\end{ack}

\bibliographystyle{amsalpha}

\end{document}